# Robust mixture regression modeling based on the Generalized M (GM)-estimation method

Fatma Zehra Doğru[1*] and Olcay Arslan[1]

[1]Ankara University, Faculty of Science, Department of Statistics, 06100 Ankara/Turkey
fzdogru@ankara.edu.tr, oarslan@ankara.edu.tr

**Abstract**

Bai (2010) and Bai et al. (2012) proposed robust mixture regression method based on the M regression estimation. However, the M-estimators are robust against the outliers in response variables, but they are not robust against the outliers in explanatory variables (leverage points). In this paper, we propose a robust mixture regression procedure to handle the outliers and the leverage points, simultaneously. Our proposed mixture regression method is based on the GM regression estimation. We give an Expectation Maximization (EM) type algorithm to compute estimates for the parameters of interest. We provide a simulation study and a real data example to assess the robustness performance of the proposed method against the outliers and the leverage points.

## 1. Introduction

Mixture regression models are widely applied in areas such as engineering, genetics, biology, econometrics and marketing, which are used to examine the relationship between variables coming from some unknown latent groups. These models were first introduced by Quandt (1972) and Quandt and Ramsey (1978) as switching regression models.

The mixture regression model can be defined as follows. Let $x$ be a $p$-dimensional vector of explanatory variables, $y$ be the response variable and $Z$ be a latent variable with $P(Z_j = i|x) = \pi_i$ for $i = 1, ..., g$, denote the mixing probabilities with $\sum_{i=1}^{g} \pi_i = 1$, $0 \leq \pi_i \leq 1$. Suppose that given $Z = i$, the response variable $y$ depends on the explanatory variable $x$ in a linear way as

$$y = x'\beta_i + \epsilon_i, \quad i = 1,2,\cdots,g, \tag{1}$$

where $\epsilon_i$ is the error term, $\beta_i = (\beta_{i1}, \beta_{i2}, \cdots, \beta_{ip})'$ is the unknown vector of regression parameters, $g$ is the number of components in mixture regression model and $x$ includes both predictors and constant 1. If it is assumed that the distributions of $\epsilon_i$'s are the member of location-scale family with zero means and $\sigma_i$ scale parameters, then the conditional density function of $y$ given $x$ is

$$f(y; x, \Theta) = \sum_{i=1}^{g} \pi_i f_i(y; x'\beta_i, \sigma_i), \tag{2}$$

where $f_i(y; x'\beta_i, \sigma_i)$ is the density function of the $ith$ component and $\Theta = (w_1, \cdots, w_g, \beta_1, \cdots, \beta_g, \sigma_1, \cdots, \sigma_g)'$ is the unknown parameter vector. This model is called as a g-component mixture regression model.

The distributions of $\epsilon_i$'s are usually assumed to be normal and under this assumption the mixture regression model given in (2) becomes a finite mixture of normal distributions. However, under normality assumption the resulting estimator can be affected by the outliers and heavy tailed error distributions. Thus, robust mixture regression procedures have been proposed to overcome these problems. For instance, Neykov et al. (2007) proposed robust fitting of mixtures based on the trimmed likelihood estimator. Markatou (2000) and Shen et al. (2004) used the weight factor for robustly estimating the parameters in a mixture regression model. Bashir and Carter (2012) used the S-estimation method to find the robust estimators for the parameters of the mixture linear regression model. Bai (2010) and Bai et al. (2012) proposed a robust estimation procedure for the mixture regression procedure based on the M regression estimation method. Further, there are several papers on the robust mixture regression model based on heavy-tailed distributions. For example, Wei (2012) and Yao et al. (2014) studied mixture regression model based on the *t* distribution, Zhang (2013) proposed the robust mixture regression model using the Pearson Type VII distribution and Song et al. (2014) explored robust mixture regression models using the mixture of Laplace distribution. Recently, Doğru (2015) and Doğru and Arslan (2015) have proposed robust mixture regression procedure based on the skew *t* distribution (Azzalini and Capitaino (2003)) to model skewness and heavy-taildness of the errors. This mixture regression model is an extension of the mixture of skew *t* distribution proposed by Lin et al. (2007a). Doğru and Arslan (2015) have also proposed a robust mixture regression procedure using the mixture of different distributions to cope with the heterogeneity in data.

Bai (2010) and Bai et al. (2012) combined the M-estimation and the EM algorithm to propose their robust mixture regression procedure. They first take the normal mixture regression model and carry out the EM algorithm using the M objective function instead of the least square criterion. By doing this they made their estimators robust against the outliers in response. However, since the M-estimators are robust against the outliers in the $y$ direction but not robust against the outliers in the $x$ direction, their estimators will be sensitive to the outliers in the $x$ direction (leverage points). The purpose of this paper is to propose a mixture regression procedure based on the GM-estimation method (Mallows (1975), Hampel (1978), Krasker (1980), Krasker and Welsch (1982), Hampel et al. (1986), Maronna et al. (2006), Jurečková and Picek (2006)) which will be robust against the outliers in the $x$ and $y$ directions simultaneously.

The paper is organized as follows. In Section 2, we give the definition of the mixture regression model based on the GM-estimation method and also give an EM type algorithm to obtain the parameter estimators. In Section 3 and 4, we give a simulation study and a real data example to compare the performance of the proposed estimation procedure over the estimation procedure proposed by Bai (2010) and Bai et al. (2012). The paper is finilazed with a conclusion section.

**2. Mixture regression model based on the GM-estimation method**

Let $\{(x_1, y_1), (x_2, y_2), \cdots, (x_n, y_n)\}$ be a sample of observations. If we assume that the error terms have the normal distribution with 0 mean and $\sigma^2$ variance in the mixture regression model, the ML estimator of $\boldsymbol{\Theta}$ for a g-component mixture regression model given in (2) will be

$$\widehat{\boldsymbol{\Theta}} = \arg\max_{\boldsymbol{\Theta}} \sum_{j=1}^{n} \log\left(\sum_{i=1}^{g} \pi_i \phi(y_j; x_j' \boldsymbol{\beta}_i, \sigma_i^2)\right). \quad (3)$$

Since there is not an explicit maximizer of (3), the EM algorithm (Dempster et al. (1977)) is usually used to obtain the ML estimates of $\boldsymbol{\Theta}$.

To carry out the EM algorithm let $z_{ij}$ be the latent variables defined as

$$z_{ij} = \begin{cases} 1, & if\ j^{th}\ observation\ is\ from\ i^{th}\ component \\ 0, & otherwise \end{cases} \quad (4)$$

where, $j = 1, \cdots, n$ and $i = 1, \cdots, g$. Since $z_{ij}$ cannot be observed they will be regarded as the missing observations and $(\mathbf{y}, \mathbf{z}_j)$ will form the complete data. Then, the complete data log-likelihood function given $\mathbf{X}$ can be written as

$$\ell_c(\boldsymbol{\Theta}|\mathbf{y}, \mathbf{z}_j) = \sum_{j=1}^{n}\sum_{i=1}^{g} z_{ij}\left(\log(\pi_i) - \frac{1}{2}\log(2\pi) - \frac{1}{2}\log(\sigma_i^2)\right) - \sum_{j=1}^{n}\sum_{i=1}^{g} z_{ij} \frac{(y_j - \mathbf{x}_j'\boldsymbol{\beta}_i)^2}{2\sigma_i^2}, \quad (5)$$

where $\mathbf{X} = (\mathbf{x}_1, \cdots, \mathbf{x}_n)'$, $\mathbf{y} = (y_1, \cdots, y_n)$ and $\mathbf{z}_j = (z_{1j}, \cdots, z_{gj})$. Since the second term of (5) is least square criterion, the resulting estimators will not be robust against in any type of outliers. To gain robustness squared term in this part can be replaced by a robust objective function such as Huber or Tukey biweight functions. This replacement will result estimators that are robust against the heavy tails and outliers in the $y$ direction. To further gain robustness against the outliers in the $x$ as well as $y$ directions we have to use a robust criterion that is used in GM-estimation method. In this paper, we will do this and chance the square term in second part of (5) with a GM criterion.

Let $\zeta(\cdot)$ be a robust criterion function. Then, we will rewrite the complete data log likelihood function for $(\mathbf{y}, \mathbf{z}_j)$ given $\mathbf{X}$ as follows

$$\ell_c(\boldsymbol{\Theta}|\mathbf{y}, \mathbf{z}_j) = \sum_{j=1}^{n}\sum_{i=1}^{g} z_{ij}\left(\log(\pi_i) - \frac{1}{2}\log(2\pi) - \frac{1}{2}\log(\sigma_i^2)\right)$$
$$- \sum_{j=1}^{n}\sum_{i=1}^{g} z_{ij}\, \sigma_i\, \zeta\left(\mathbf{x}_j, \frac{y_j - \mathbf{x}_j'\boldsymbol{\beta}_i}{\sigma_i}\right). \quad (6)$$

This function instead of function given in (5) will be used to develop the EM algorithm to get estimates for the parameters. Note that the function $\zeta\left(\mathbf{x}, \frac{r}{\sigma}\right), r = y - \mathbf{x}'\boldsymbol{\beta}$, depends on $\mathbf{x}$ and standardized residuals. In literature, there are two forms of $\zeta(\cdot)$ used to obtain the GM- estimators. If $\zeta(\mathbf{x}, r) = w(\mathbf{x})\rho(r)$, the estimator is called Mallows type (Mallows (1975)) and if $\zeta(\mathbf{x}, r) = w^2(\mathbf{x})\rho\left(\frac{r}{w(\mathbf{x})}\right)$, the estimator is called Schweppe type (Handschin et al. (1975), Hampel et al. (1986)). Here, $w(\mathbf{x})$ is a decreasing weight function of Mahalanobis distance of $\mathbf{x}$ and it is designed to reduce the effect of leverage points on the resulting estimators. If we take $w(\mathbf{x}) = 1$, we get the M-estimators. There are many different weight functions used in literature. In this study, we will use the following weight function proposed by Simpson et al. (1992)

$$w(\mathbf{x}_j) = \min\left(1, \left(\frac{b}{d_n(\mathbf{x}_j)}\right)^{1/2}\right), \quad (7)$$

where $d_n(\mathbf{x}_j) = (\mathbf{x}_j - m(\mathbf{X}))' C(\mathbf{X})^{-1} (\mathbf{x}_j - m(\mathbf{X}))$ is the robust Mahalanobis distance based on minumum covariance determinant (MCD) (Rousseeuw and Van Driessen (1999)) estimators of sample mean and variance covariance matrix of $\mathbf{x}_j$ and $b$ is the $(1 - \gamma)$ quantile of the chi-squared distribution with $(p - 1)$ degrees of freedom.

Now, we give the steps of the EM type algorithm for the GM based mixture regression model. First we have to take the conditional expectation of the complete data log-likelihood function to get rid of the latency of $z_{ij}$

$$E(\ell_c(\mathbf{\Theta}; \mathbf{y}, \mathbf{z}_j)|y_j) = \sum_{j=1}^{n}\sum_{i=1}^{g} E(z_{ij}|y_j)\left(\log(\pi_i) - \frac{1}{2}\log 2\pi - \frac{1}{2}\log \sigma_i^2\right)$$
$$- \sum_{j=1}^{n}\sum_{i=1}^{g} E(z_{ij}|y_j)\sigma_i \zeta\left(\mathbf{x}_j, \frac{y_j - \mathbf{x}_j'\boldsymbol{\beta}_i}{\sigma_i}\right). \tag{8}$$

If $\zeta$ is differentiable, taking the derivative of (8) with respect to $\boldsymbol{\beta}_i$ and setting to zero gives the following equation

$$\sum_{j=1}^{n} E(z_{ij}|y_j)\mathbf{x}_j \eta\left(\mathbf{x}_j, \frac{y_j - \mathbf{x}_j'\boldsymbol{\beta}_i}{\sigma_i}\right) = 0, \tag{9}$$

where $\eta = \zeta'$ and $\eta(\mathbf{x}, r) = w(\mathbf{x})\psi(r)$ for Mallows type estimator and $\eta(\mathbf{x}, r) = w(\mathbf{x})\psi\left(\frac{r}{w(\mathbf{x})}\right)$ for Schweppe type estimator. Note that $\psi$ is an appropriate function defined on real number. Solving this equation will give the estimator for $\boldsymbol{\beta}_i$. However, since this estimator will be a function of $E(z_{ij}|y_j)$ it cannot be used. To make it functional we have to compute this conditional expectation given $y_j$ and the current estimate $\widehat{\mathbf{\Theta}}$, which will be

$$\hat{z}_{ij} = E(z_{ij}|y_j, \widehat{\mathbf{\Theta}}) = \frac{\hat{\pi}_i \phi(y_j; \mathbf{x}_j'\widehat{\boldsymbol{\beta}}_i, \hat{\sigma}_i^2)}{\sum_{i=1}^{g} \hat{\pi}_i \phi(y_j; \mathbf{x}_j'\widehat{\boldsymbol{\beta}}_i, \hat{\sigma}_i^2)}. \tag{10}$$

Then after some straightforward algebra, we obtain

$$\sum_{j=1}^{n} \hat{z}_{ij}\, \mathbf{x}_j \eta\left(\mathbf{x}_j, \frac{y_j - \mathbf{x}_j'\boldsymbol{\beta}_i}{\sigma_i}\right) = \sum_{j=1}^{n} \hat{z}_{ij}\, w(\mathbf{x}_j)\mathbf{x}_j W^*\left(\frac{y_j - \mathbf{x}_j'\widehat{\boldsymbol{\beta}}_i}{\hat{\sigma}_i}\right)\left(\frac{y_j - \mathbf{x}_j'\boldsymbol{\beta}_i}{\sigma_i}\right)$$
$$= \sum_{j=1}^{n} \hat{z}_{ij}^*\, w(\mathbf{x}_j)\mathbf{x}_j \left(\frac{y_j - \mathbf{x}_j'\boldsymbol{\beta}_i}{\sigma_i}\right), \tag{11}$$

where

$$\hat{z}_{ij}^* = \hat{z}_{ij} W^*\left(\frac{y_j - \mathbf{x}_j'\widehat{\boldsymbol{\beta}}_i}{\hat{\sigma}_i}\right)$$

and $W^*\left(\frac{y_j - \mathbf{x}_j'\widehat{\boldsymbol{\beta}}_i}{\hat{\sigma}_i}\right)$ shows the weight function for Mallows and Schweppe type GM-estimation methods. Weights for Mallows estimator are $\psi\left(\frac{y_j - \mathbf{x}_j'\widehat{\boldsymbol{\beta}}_i}{\hat{\sigma}_i}\right)\Big/\left(\frac{y_j - \mathbf{x}_j'\widehat{\boldsymbol{\beta}}_i}{\hat{\sigma}_i}\right)$ and for Schweppe type estimator are $\psi\left(\frac{y_j - \mathbf{x}_j'\widehat{\boldsymbol{\beta}}_i}{\hat{\sigma}_i w(\mathbf{x}_j)}\right)\Big/\left(\frac{y_j - \mathbf{x}_j'\widehat{\boldsymbol{\beta}}_i}{\hat{\sigma}_i}\right)$. Here $\sigma_i$ will be estimated using the M-scale estimation method (Huber and Ronchetti (2009), Maronna et al. (2006), Jurečková and Picek (2006)) which is adapted for the mixture regression model.

Solving (11) with respect to $\boldsymbol{\beta}_i$ will give

$$\widehat{\boldsymbol{\beta}}_i = \left(\sum_{j=1}^{n} \hat{z}_{ij}^* w(\mathbf{x}_j)\mathbf{x}_j \mathbf{x}_j'\right)^{-1}\left(\sum_{j=1}^{n} \hat{z}_{ij}^* w(\mathbf{x}_j)\mathbf{x}_j y_j\right). \tag{12}$$

Now, the steps of EM type algorithm can be given as follows:

**EM algorithm:**

**1.** Take initial parameter estimates $\boldsymbol{\Theta}^{(0)}$ and fix a stopping rule $\Delta$.
**2. E-step:** Given the current parameter values $\widehat{\boldsymbol{\Theta}}^{(k)}$, compute the following conditional expectation for $k = 0,1,2,\cdots$

$$\hat{z}_{ij}^{(k)} = E(Z_{ij}|y_j, \widehat{\boldsymbol{\Theta}}^{(k)}) = \frac{\hat{\pi}_i^{(k)} \phi\left(y_j; \boldsymbol{x}_j'\widehat{\boldsymbol{\beta}}_i^{(k)}, \hat{\sigma}_i^{2(k)}\right)}{\sum_{i=1}^{g} \hat{\pi}_i^{(k)} \phi\left(y_j; \boldsymbol{x}_j'\widehat{\boldsymbol{\beta}}_i^{(k)}, \hat{\sigma}_i^{2(k)}\right)} \quad . \tag{13}$$

**3. M-step:**
i) Update the mixing probabilities using

$$\hat{\pi}_i^{(k+1)} = \frac{\sum_{j=1}^{n} \hat{z}_{ij}^{(k)}}{n} . \tag{14}$$

ii) Using equation (12) compute $\widehat{\boldsymbol{\beta}}_i^{(k+1)}$

$$\widehat{\boldsymbol{\beta}}_i^{(k+1)} = \left(\sum_{j=1}^{n} \hat{z}_{ij}^{*(k)} w(\boldsymbol{x}_j) \boldsymbol{x}_j \boldsymbol{x}_j'\right)^{-1} \left(\sum_{j=1}^{n} \hat{z}_{ij}^{*(k)} w(\boldsymbol{x}_j) \boldsymbol{x}_j y_j\right). \tag{15}$$

iii) Calculate $\hat{\sigma}_i^{2(k+1)}$ using

$$\hat{\sigma}_i^{2(k+1)} = \frac{\hat{\sigma}_i^{2(k)}}{a \sum_{j=1}^{n} \hat{z}_{ij}^{(k)}} \sum_{j=1}^{n} \hat{z}_{ij}^{(k)} \chi\left(\frac{y_j - \boldsymbol{x}_j'\widehat{\boldsymbol{\beta}}_i^{(k)}}{\hat{\sigma}_i^{(k)}}\right), \tag{16}$$

where

$$\chi\left(\frac{y_j - \boldsymbol{x}_j'\widehat{\boldsymbol{\beta}}_i^{(k)}}{\hat{\sigma}_i^{(k)}}\right) = \psi\left(\frac{y_j - \boldsymbol{x}_j'\widehat{\boldsymbol{\beta}}_i^{(k)}}{\hat{\sigma}_i^{(k)}}\right) \frac{\left(y_j - \boldsymbol{x}_j'\widehat{\boldsymbol{\beta}}_i^{(k)}\right)}{\hat{\sigma}_i^{(k)}} - \rho\left(\frac{y_j - \boldsymbol{x}_j'\widehat{\boldsymbol{\beta}}_i^{(k)}}{\hat{\sigma}_i^{(k)}}\right),$$

$$a = \frac{n-p}{n} E_\Phi\left(\chi\left(\frac{y_j - \boldsymbol{x}_j'\widehat{\boldsymbol{\beta}}_i^{(k)}}{\hat{\sigma}_i^{(k)}}\right)\right).$$

**4.** Repeat E and M steps until the convergence criteria $\|\boldsymbol{\Theta}^{(k+1)} - \boldsymbol{\Theta}^{(k)}\| < \Delta$ is satisfied.

## 3. Simulation study

In this section, we give a simulation study to compare the performance of the proposed estimators (MixregGM-Mallows and MixregGM-Schweppe) with the mixture regression estimators based on M-estimation (Mixreg-Huber and Mixreg-Tukey) proposed by Bai (2010) and Bai et al. (2012) in terms of bias and mean square error (MSE). The bias and MSE are computed using the following formulas

$$\widehat{bias}(\hat{\theta}) = \bar{\theta} - \theta , \widehat{MSE}(\hat{\theta}) = \frac{1}{N}\sum_{i=1}^{N}(\hat{\theta}_i - \theta)^2$$

where $\theta$ is the true parameter value, $\hat{\theta}_i$ is the estimation of $\theta$ from the ith simulated data and $\bar{\theta} = \frac{1}{N}\sum_{i=1}^{N}\hat{\theta}_i$. The number of replications ($N$) is 500. We take sample sizes as 200 and 400 for all simulation settings.

For the simulation, we use Huber's $\psi$ function $\psi_c(x) = \max(-c, \min(c, x))$ with $c = 1.345$ and we choose $\eta(x, r) = w(x)\psi(r)$ for Mallows type estimator and $\eta(x, r) = w(x)\psi\left(\frac{r}{w(x)}\right)$ for Schweppe type estimator. The weights $w(x)$ are computed using (7). For the Tukey's bisquare function $\psi_c(x) = x(1 - (x/c)^2)^2 I(|x| \leq c)$ which will be used to compute the mixture regression M estimator proposed by Bai (2010) and Bai et al. (2012), the tuning constant will be taken as 4.685. All simulation studies are conducted using MATLAB R2013a.

We consider two simulation scenarios.

**Scenario 1.** We generate the data $\{(x_j, y_j), j = 1, ..., n\}$ from a two-component mixture regression model (Bai (2010))

$$y = \begin{cases} 0 + 4x + \epsilon_1, & z = 1, \\ 0 - 4x + \epsilon_2, & z = 2, \end{cases}$$

where $P(z = 1) = 0.5 = \pi_1$ and $x \sim N(0,1)$. Furthermore, the model coefficients are $\boldsymbol{\beta}_1 = (\beta_{10}, \beta_{11})' = (0,4)'$ and $\boldsymbol{\beta}_2 = (\beta_{20}, \beta_{21})' = (0,-4)'$. For the error distribution, we take the following cases:

Case I: $\epsilon_1, \epsilon_2 \sim N(0,1)$, standard normal distribution.
Case II: $\epsilon_1, \epsilon_2 \sim t_4$, $t$ distribution with the degrees of freedom 4.
Case III: $\epsilon_1, \epsilon_2 \sim 0.95N(0,1) + 0.05N(0,25)$, contaminated normal distribution.
Case IV: $\epsilon_1, \epsilon_2 \sim N(0,1)$, standard normal distribution with outliers, (we add 5 and 10 outliers in the $x$ direction for $n = 200$ and 400, respectively).

Tables 1 and 2 display the simulation results for the Scenario 1. In the tables we give bias and MSE values. We observe the followings from the simulation results. Although, for the Cases I, II and III all the estimators have similar behavior, MixregGM-Schweppe has smaller MSE values in most simulation conditions. For the Case IV, which is the case contains outliers (5 outliers for $n = 200$ and 10 outliers for $n = 400$) in the $x$ direction, we observe that Mixreg-Huber and Mixreg-Tukey estimators are drastically affected by the outliers. On the other hand, MixregGM-Mallows and MixregGM-Schweppe estimators have the lowest bias and MSE values in almost all simulation conditions.

**Table 1.** MSE (bias) values of estimates for $n = 200$.

|  | Mixreg-Huber | Mixreg-Tukey | MixregGM-Mallows | MixregGM-Schweppe |
|---|---|---|---|---|
| | Case I: $\epsilon_1, \epsilon_2 \sim N(0,1)$ | | | |
| $\beta_{10}: 0$ | 0.0119 (0.0008) | 0.0113 (0.0012) | 0.0118 (0.0019) | 0.0118 (0.0013) |
| $\beta_{20}: 0$ | 0.0126 (0.0024) | 0.0122 (0.0030) | 0.0127 (0.0017) | 0.0127 (0.0019) |
| $\beta_{11}: 4$ | 0.0113 (-0.0101) | 0.0108 (-0.0094) | 0.0114 (-0.0097) | 0.0115 (-0.0099) |
| $\beta_{21}: -4$ | 0.0120 (0.0058) | 0.0117 (0.0074) | 0.0120 (0.0065) | 0.0120 (0.0057) |
| $\pi_1: 0.5$ | 0.0016 (0.0020) | 0.0016 (0.0019) | 0.0016 (0.0017) | 0.0016 (0.0017) |
| | Case II: $\epsilon_1, \epsilon_2 \sim t_4$ | | | |
| $\beta_{10}: 0$ | 0.0182 (0.0027) | 0.0211 (0.0006) | 0.0185 (0.0039) | 0.0184 (0.0042) |
| $\beta_{20}: 0$ | 0.0181 (0.0065) | 0.0251 (0.0073) | 0.0192 (0.0049) | 0.0191 (0.0047) |
| $\beta_{11}: 4$ | 0.0159 (0.0089) | 0.0185 (0.0175) | 0.0167 (0.0021) | 0.0165 (0.0018) |
| $\beta_{21}: -4$ | 0.0145 (-0.0022) | 0.0188 (-0.0090) | 0.0169 (0.0084) | 0.0165 (0.0084) |
| $\pi_1: 0.5$ | 0.0016 (-0.0021) | 0.0016 (-0.0021) | 0.0016 (-0.0025) | 0.0016 (-0.0025) |
| | Case III: $\epsilon_1, \epsilon_2 \sim 0.95N(0,1) + 0.05N(0,25)$ | | | |
| $\beta_{10}: 0$ | 0.0135 (0.0026) | 0.0152 (-0.0001) | 0.0150 (0.0033) | 0.0149 (0.0035) |
| $\beta_{20}: 0$ | 0.0161 (0.0050) | 0.0176 (0.0055) | 0.0170 (0.0045) | 0.0170 (0.0044) |
| $\beta_{11}: 4$ | 0.0124 (0.0101) | 0.0139 (0.0140) | 0.0134 (-0.0009) | 0.0130 (-0.0008) |
| $\beta_{21}: -4$ | 0.0132 (-0.0105) | 0.0145 (-0.0127) | 0.0146 (0.0036) | 0.0143 (0.0033) |
| $\pi_1: 0.5$ | 0.0015 (-0.0005) | 0.0015 (-0.0005) | 0.0018 (-0.0008) | 0.0018 (-0.0008) |
| | Case IV: $\epsilon_1, \epsilon_2 \sim N(0,1)$ (5 outliers) | | | |
| $\beta_{10}: 0$ | 0.2043 (-0.1735) | 0.1860 (-0.1685) | 0.0643 (-0.0681) | 0.0628 (-0.0686) |
| $\beta_{20}: 0$ | 0.0194 (0.0025) | 0.0179 (0.0025) | 0.0158 (0.0015) | 0.0157 (0.0018) |
| $\beta_{11}: 4$ | 13.5129 (-3.6757) | 13.3367 (-3.6517) | 3.1684 (-1.7710) | 2.8973 (-1.6914) |
| $\beta_{21}: -4$ | 0.0172 (-0.0285) | 0.0162 (-0.0291) | 0.0171 (0.0403) | 0.0172 (0.0423) |
| $\pi_1: 0.5$ | 0.0123 (0.0964) | 0.0117 (0.0938) | 0.0020 (-0.0071) | 0.0021 (-0.0116) |

Note: Value in parentheses indicates the bias

**Table 2.** MSE (bias) values of estimates for $n = 400$.

|  | Mixreg-Huber | Mixreg-Tukey | MixregGM-Mallows | MixregGM-Schweppe |
|---|---|---|---|---|
| | Case I: $\epsilon_1, \epsilon_2 \sim N(0,1)$ | | | |
| $\beta_{10}: 0$ | 0.0058 (0.0031) | 0.0054 (0.0035) | 0.0057 (0.0031) | 0.0057 (0.0031) |
| $\beta_{20}: 0$ | 0.0063 (-0.0009) | 0.0059 (-0.0015) | 0.0063 (-0.0010) | 0.0063 (-0.0009) |
| $\beta_{11}: 4$ | 0.0053 (-0.0023) | 0.0050 (-0.0031) | 0.0054 (-0.0019) | 0.0055 (-0.0021) |
| $\beta_{21}: -4$ | 0.0059 (-0.0013) | 0.0056 (-0.0007) | 0.0059 (-0.0012) | 0.0059 (-0.0013) |
| $\pi_1: 0.5$ | 0.0008 (0.0012) | 0.0008 (0.0012) | 0.0008 (0.0012) | 0.0008 (0.0012) |
| | Case II: $\epsilon_1, \epsilon_2 \sim t_4$ | | | |
| $\beta_{10}: 0$ | 0.0090 (0.0029) | 0.0120 (0.0063) | 0.0094 (0.0042) | 0.0094 (0.0042) |
| $\beta_{20}: 0$ | 0.0086 (-0.0053) | 0.0107 (-0.0051) | 0.0089 (-0.0056) | 0.0089 (-0.0055) |
| $\beta_{11}: 4$ | 0.0066 (0.0138) | 0.0094 (0.0222) | 0.0070 (0.0063) | 0.0068 (0.0062) |
| $\beta_{21}: -4$ | 0.0073 (-0.0036) | 0.0084 (-0.0109) | 0.0076 (0.0042) | 0.0075 (0.0044) |
| $\pi_1: 0.5$ | 0.0010 (-0.0002) | 0.0010 (-0.0002) | 0.0010 (-0.0002) | 0.0010 (-0.0002) |
| | Case III: $\epsilon_1, \epsilon_2 \sim 0.95N(0,1) + 0.05N(0,25)$ | | | |
| $\beta_{10}: 0$ | 0.0072 (0.0002) | 0.0085 (0.0005) | 0.0075 (0.0008) | 0.0074 (0.0008) |
| $\beta_{20}: 0$ | 0.0077 (-0.0077) | 0.0099 (-0.0077) | 0.0079 (-0.0087) | 0.0078 (-0.0087) |
| $\beta_{11}: 4$ | 0.0066 (0.0109) | 0.0072 (0.0133) | 0.0066 (-0.0027) | 0.0065 (-0.0020) |
| $\beta_{21}: -4$ | 0.0064 (-0.0095) | 0.0081 (-0.0157) | 0.0072 (0.0015) | 0.0070 (0.0015) |
| $\pi_1: 0.5$ | 0.0008 (-0.0006) | 0.0008 (-0.0005) | 0.0009 (-0.0003) | 0.0009 (-0.0003) |
| | Case IV: $\epsilon_1, \epsilon_2 \sim N(0,1)$ (10 outliers) | | | |
| $\beta_{10}: 0$ | 0.1067 (-0.1820) | 0.1024 (-0.1802) | 0.0330 (-0.0777) | 0.0331 (-0.0791) |
| $\beta_{20}: 0$ | 0.0089 (0.0051) | 0.0083 (0.0052) | 0.0082 (0.0018) | 0.0083 (0.0017) |
| $\beta_{11}: 4$ | 13.0840 (-3.6171) | 12.9210 (-3.5944) | 3.1303 (-1.7644) | 2.8640 (-1.6868) |
| $\beta_{21}: -4$ | 0.0080 (-0.0276) | 0.0077 (-0.0286) | 0.0081 (0.0368) | 0.0084 (0.0388) |
| $\pi_1: 0.5$ | 0.0090 (0.0878) | 0.0087 (0.0863) | 0.0011 (-0.0099) | 0.0012 (-0.0143) |

Note: Value in parentheses indicates the bias

**Scenario 2.** We generate the data $\{(x_{1j}, x_{2j}, y_j), j = 1, \ldots, n\}$ from the following two component mixture regression models (Bai et al. (2012))

$$y = \begin{cases} 0 + x_1 + x_2 + \epsilon_1, & z = 1, \\ 0 - x_1 - x_2 + \epsilon_2, & z = 2, \end{cases}$$

where $P(z = 1) = 0.25 = \pi_1$, $x_1 \sim N(0,1)$ and $x_2 \sim N(0,1)$. Here, the model coefficients are $\boldsymbol{\beta}_1 = (\beta_{10}, \beta_{11}, \beta_{12})' = (0,1,1)'$ and $\boldsymbol{\beta}_2 = (\beta_{20}, \beta_{21}, \beta_{22})' = (0,-1,-1)'$. We explore the following error distributions:

Case I: $\epsilon_1, \epsilon_2 \sim N(0,1)$, standard normal distribution.
Case II: $\epsilon_1, \epsilon_2 \sim t_3$, $t$ distribution with the degrees of freedom 3.
Case III: $\epsilon_1, \epsilon_2 \sim 0.95N(0,1) + 0.05N(0,25)$, contaminated normal distribution.
Case IV: $\epsilon_1, \epsilon_2 \sim N(0,1)$, standard normal distribution with outliers, (we add 5 outliers for $n = 200$ and 10 outliers for $n = 400$ in the $x$ direction).

In Tables 3 and 4, we give the simulation results. The tables show the bias and MSE values of the estimates. The simulation results show that all estimators have similar performance when the error terms have normal, heavy-tailed and contaminated normal distributions. However, when we add outliers (5 outliers for $n = 200$ and 10 outliers for $n = 400$) to the data in the $x$ direction and error terms have the normal distribution, Mixreg-Huber and Mixreg-Tukey estimators are failed to find the right groups as they are drastically influenced by the outliers. On the contrary, since the MixregGM-Mallows and MixregGM-Schweppe estimators are resistant to the outliers in the $x$ direction, they have the lowest bias and MSE values for almost all the estimators compare to the others.

In summary, after evaluating the results of all simulation studies, we observe that mixture regression model based on the robust estimation method should be used when the data includes outliers. In particular, if the data have the outliers in the $x$ direction, since the mixture regression model based on the GM-estimation method gives reliable results, the mixture regression GM-estimators should be used instead of the mixture regression model based on the M-estimation proposed by Bai (2010) and Bai et al. (2012).

**Table 3.** MSE (bias) values of estimates for $n = 200$.

|  | Mixreg-Huber | Mixreg-Tukey | MixregGM-Mallows | MixregGM-Schweppe |
|---|---|---|---|---|
| Case I: $\epsilon_1, \epsilon_2 \sim N(0,1)$ | | | | |
| $\beta_{10}: 0$ | 0.0556 (-0.0136) | 0.0523 (-0.0122) | 0.0533 (-0.0127) | 0.0532 (-0.0131) |
| $\beta_{20}: 0$ | 0.0098 (0.0033) | 0.0092 (0.0028) | 0.0096 (0.0034) | 0.0095 (0.0037) |
| $\beta_{11}: 1$ | 0.0512 (0.0218) | 0.0553 (0.0181) | 0.0448 (0.0265) | 0.0452 (0.0271) |
| $\beta_{21}: -1$ | 0.0095 (-0.0060) | 0.0092 (-0.0069) | 0.0087 (-0.0055) | 0.0087 (-0.0051) |
| $\beta_{12}: 1$ | 0.0534 (-0.0023) | 0.0558 (-0.0066) | 0.0473 (0.0007) | 0.0478 (0.0010) |
| $\beta_{22}: -1$ | 0.0098 (0.0031) | 0.0097 (0.0011) | 0.0095 (0.0067) | 0.0095 (0.0067) |
| $\pi_1: 0.25$ | 0.0039 (0.0024) | 0.0041 (0.0046) | 0.0028 (0.0011) | 0.0028 (0.0010) |
| Case II: $\epsilon_1, \epsilon_2 \sim t_3$ | | | | |
| $\beta_{10}: 0$ | 0.0940 (-0.0111) | 0.1061 (-0.0222) | 0.1307 (-0.0207) | 0.1316 (-0.0179) |
| $\beta_{20}: 0$ | 0.0208 (-0.0007) | 0.0213 (0.0028) | 0.0223 (-0.0017) | 0.0167 (-0.0027) |
| $\beta_{11}: 1$ | 0.3679 (-0.2411) | 0.3694 (-0.2376) | 0.2932 (-0.2318) | 0.2491 (-0.1867) |
| $\beta_{21}: -1$ | 0.0238 (-0.0450) | 0.0254 (-0.0596) | 0.0329 (0.0055) | 0.0174 (-0.0089) |
| $\beta_{12}: 1$ | 0.3597 (-0.2580) | 0.3648 (-0.2526) | 0.3032 (-0.2741) | 0.2476 (-0.2372) |
| $\beta_{22}: -1$ | 0.0247 (-0.0466) | 0.0280 (-0.0608) | 0.0463 (0.0119) | 0.0213 (-0.0048) |
| $\pi_1: 0.25$ | 0.0333 (0.1002) | 0.0290 (0.0971) | 0.0130 (0.0487) | 0.0104 (0.0429) |
| Case III: $\epsilon_1, \epsilon_2 \sim 0.95N(0,1) + 0.05N(0,25)$ | | | | |
| $\beta_{10}: 0$ | 0.0681 (0.0088) | 0.0857 (0.0107) | 0.0942 (0.0080) | 0.0857 (0.0026) |
| $\beta_{20}: 0$ | 0.0107 (-0.0008) | 0.0113 (-0.0002) | 0.0289 (0.0058) | 0.0292 (0.0060) |
| $\beta_{11}: 1$ | 0.0918 (0.0054) | 0.1102 (0.0099) | 0.1854 (-0.1174) | 0.1467 (-0.1218) |
| $\beta_{21}: -1$ | 0.0122 (-0.0064) | 0.0125 (-0.0141) | 0.0235 (0.0168) | 0.0229 (0.0154) |
| $\beta_{12}: 1$ | 0.1088 (-0.0269) | 0.1483 (-0.0460) | 0.1963 (-0.1823) | 0.1842 (-0.1746) |
| $\beta_{22}: -1$ | 0.0116 (-0.0088) | 0.0118 (-0.0117) | 0.0152 (0.0158) | 0.0149 (0.0140) |
| $\pi_1: 0.25$ | 0.0061 (0.0089) | 0.0058 (0.0126) | 0.0063 (0.0106) | 0.0062 (0.0114) |
| Case IV: $\epsilon_1, \epsilon_2 \sim N(0,1)$ (5 outliers) | | | | |
| $\beta_{10}: 0$ | 0.3503 (-0.0404) | 0.2706 (-0.0600) | 0.0611 (-0.0179) | 0.0612 (-0.0136) |
| $\beta_{20}: 0$ | 0.0971 (-0.0022) | 0.1250 (0.0030) | 0.0128 (0.0003) | 0.0129 (0.0007) |
| $\beta_{11}: 1$ | 2.1562 (-1.2028) | 2.3962 (-1.3143) | 0.5132 (-0.5245) | 0.4968 (-0.4955) |
| $\beta_{21}: -1$ | 0.8802 (0.7992) | 0.8348 (0.7484) | 0.0145 (-0.0446) | 0.0145 (-0.0434) |
| $\beta_{12}: 1$ | 2.5613 (-1.4266) | 2.7394 (-1.5047) | 0.5371 (-0.5470) | 0.5244 (-0.5190) |
| $\beta_{22}: -1$ | 1.0061 (0.8581) | 0.9528 (0.8031) | 0.0167 (-0.0528) | 0.0166 (-0.0506) |
| $\pi_1: 0.25$ | 0.0707 (0.1609) | 0.0777 (0.2008) | 0.0228 (0.1154) | 0.0222 (0.1100) |

Note: Value in parentheses indicates the bias

**Table 4.** MSE (bias) values of estimates for $n = 400$.

|  | Mixreg-Huber | Mixreg-Tukey | MixregGM-Mallows | MixregGM-Schweppe |
|---|---|---|---|---|
| | Case I: $\epsilon_1, \epsilon_2 \sim N(0,1)$ | | | |
| $\beta_{10}: 0$ | 0.0205 (-0.0048) | 0.0196 (-0.0047) | 0.0205 (-0.0030) | 0.0206 (-0.0030) |
| $\beta_{20}: 0$ | 0.0043 (-0.0006) | 0.0041 (-0.0013) | 0.0043 (-0.0006) | 0.0043 (-0.0006) |
| $\beta_{11}: 1$ | 0.0218 (0.0034) | 0.0213 (0.0022) | 0.0198 (0.0054) | 0.0197 (0.0054) |
| $\beta_{21}: -1$ | 0.0045 (0.0051) | 0.0045 (0.0051) | 0.0046 (0.0055) | 0.0047 (0.0059) |
| $\beta_{12}: 1$ | 0.0215 (0.0057) | 0.0214 (0.0032) | 0.0201 (0.0067) | 0.0201 (0.0067) |
| $\beta_{22}: -1$ | 0.0045 (-0.0020) | 0.0044 (-0.0023) | 0.0044 (-0.0012) | 0.0044 (-0.0011) |
| $\pi_1: 0.25$ | 0.0014 (0.0003) | 0.0014 (0.0007) | 0.0012 (-0.0007) | 0.0012 (-0.0007) |
| | Case II: $\epsilon_1, \epsilon_2 \sim t_3$ | | | |
| $\beta_{10}: 0$ | 0.0415 (0.0177) | 0.0489 (0.0151) | 0.0523 (0.0129) | 0.0504 (0.0122) |
| $\beta_{20}: 0$ | 0.0080 (-0.0042) | 0.0090 (-0.0045) | 0.0071 (-0.0008) | 0.0071 (-0.0011) |
| $\beta_{11}: 1$ | 0.2600 (-0.2071) | 0.2941 (-0.2218) | 0.2283 (-0.3035) | 0.2086 (-0.2734) |
| $\beta_{21}: -1$ | 0.0105 (-0.0447) | 0.0147 (-0.0602) | 0.0090 (-0.0006) | 0.0084 (-0.0044) |
| $\beta_{12}: 1$ | 0.2644 (-0.2324) | 0.3023 (-0.2507) | 0.2456 (-0.3275) | 0.2168 (-0.2911) |
| $\beta_{22}: -1$ | 0.0104 (-0.0479) | 0.0131 (-0.0621) | 0.0115 (-0.0004) | 0.0111 (-0.0053) |
| $\pi_1: 0.25$ | 0.0255 (0.0952) | 0.0255 (0.0977) | 0.0094 (0.0613) | 0.0090 (0.0595) |
| | Case III: $\epsilon_1, \epsilon_2 \sim 0.95N(0,1) + 0.05N(0,25)$ | | | |
| $\beta_{10}: 0$ | 0.0301 (-0.0090) | 0.0360 (-0.0077) | 0.0440 (-0.0035) | 0.0439 (-0.0041) |
| $\beta_{20}: 0$ | 0.0057 (-0.0018) | 0.0060 (-0.0021) | 0.0058 (-0.0020) | 0.0059 (-0.0018) |
| $\beta_{11}: 1$ | 0.0278 (-0.0030) | 0.0356 (0.0030) | 0.1202 (-0.1988) | 0.1077 (-0.1866) |
| $\beta_{21}: -1$ | 0.0052 (-0.0075) | 0.0055 (-0.0124) | 0.0097 (0.0162) | 0.0059 (0.0136) |
| $\beta_{12}: 1$ | 0.0311 (0.0180) | 0.0423 (0.0216) | 0.1293 (-0.2092) | 0.1125 (-0.1947) |
| $\beta_{22}: -1$ | 0.0051 (-0.0125) | 0.0054 (-0.0170) | 0.0094 (0.0130) | 0.0059 (0.0101) |
| $\pi_1: 0.25$ | 0.0019 (0.0045) | 0.0019 (0.0058) | 0.0029 (0.0174) | 0.0025 (0.0154) |
| | Case IV: $\epsilon_1, \epsilon_2 \sim N(0,1)$ (10 outliers) | | | |
| $\beta_{10}: 0$ | 0.3835 (0.0181) | 0.2072 (-0.0117) | 0.0279 (-0.0122) | 0.0282 (-0.0118) |
| $\beta_{20}: 0$ | 0.1146 (-0.0077) | 0.1142 (0.0045) | 0.0069 (-0.0012) | 0.0070 (-0.0014) |
| $\beta_{11}: 1$ | 2.2019 (-1.2476) | 2.6361 (-1.4434) | 0.5449 (-0.5805) | 0.5254 (-0.5538) |
| $\beta_{21}: -1$ | 0.9721 (0.9269) | 0.9231 (0.8478) | 0.0094 (-0.0457) | 0.0095 (-0.0460) |
| $\beta_{12}: 1$ | 2.6391 (-1.4760) | 2.8753 (-1.5700) | 0.5586 (-0.5961) | 0.5487 (-0.5726) |
| $\beta_{22}: -1$ | 1.0954 (0.9869) | 0.9645 (0.8701) | 0.0088 (-0.0443) | 0.0089 (-0.0432) |
| $\pi_1: 0.25$ | 0.0671 (0.1334) | 0.0768 (0.1957) | 0.0238 (0.1251) | 0.0235 (0.1209) |

Note: Value in parentheses indicates the bias

## 4. Real data example

In this section, we will analyze the ethanol data set which is given by Hurvich et al. (1998). This data set contains the concentration of nitric oxide in engine exhaust and equivalence ratio, which is the richness of the air-ethanol mix in an engine. This data set can be accessed by using *locfit* package (Loader (1999)) in R. This data set also used by Hurn et al. (2003) in the context of mixture regression model for the Bayesian approach. Figure 1 shows the scatter plot and histogram of equivalence ratio. From these plots it is clear that there are two separate groups in the data set.

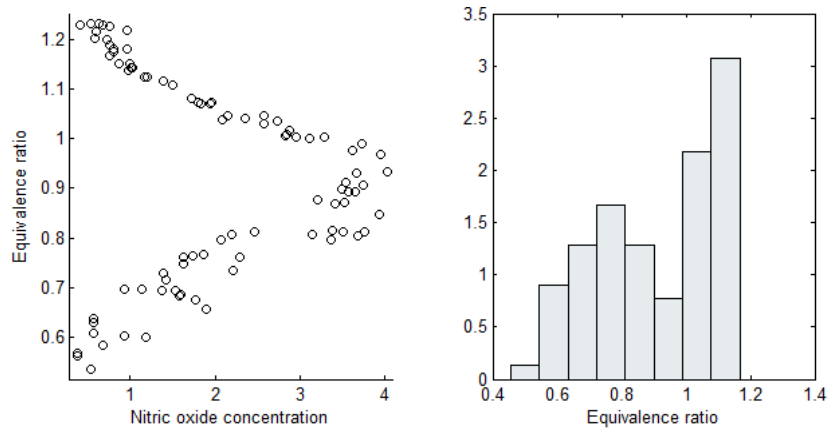

**Figure 1.** (a) The scatter plot of the data. (b) Histogram of the equivalence ratio

Using this data set, we compare the performances of the estimators in the case of with and without outliers. We show the scatter plots with the fitted regression lines obtained from Mixreg-Huber, Mixreg-Tukey, MixregGM-Mallows and MixregGM-Schweppe procedures in Figure 2. Also, we give the estimates for regression coefficients and mixing probability along with the standard errors of estimators in Table 5. The standard errors for the estimates are computed using the asymptotic covariance matrix of the estimators given in Appendix. For the Mixreg-Huber and Mixreg-Tukey, we use the asymptotic results given in Bai et al. (2012). We also include the values of the integrated complete likelihood (ICL) (Biernacki et al. (2000)) criterion for these procedures. We see from figure that all estimation procedures give similar fits, but MixregGM-Mallows gives the best result in terms of ICL criterion. We can also compute the estimates for $\sigma$ using the equation given in (16). For simplicity we fix the $\sigma$ to compute the standard errors of regression coefficients and mixing probabilities.

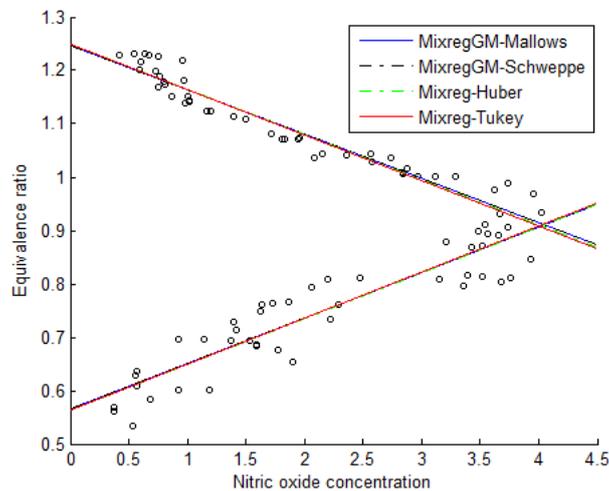

**Figure 2.** Fitted mixture regression lines for the ethanol data set

**Table 5.** Parameter estimates, standard errors of estimates and ICL information criterion for fitting different mixture regression models to the ethanol data set

| Model | Parameter Estimates | | | | | | | Information Criterion | |
|---|---|---|---|---|---|---|---|---|---|
| | $\hat{\pi}_1$ | $\hat{\beta}_{10}$ | $\hat{\beta}_{20}$ | $\hat{\beta}_{11}$ | $\hat{\beta}_{21}$ | $\hat{\sigma}_1$ | $\hat{\sigma}_2$ | $\ell_c(\hat{\Theta})$ | ICL |
| Mixreg-Huber | 0.48761 (0.00786) | 0.56468 (0.12424) | 0.08533 (0.05547) | 1.24774 (0.11663) | -0.08443 (0.05601) | 0.05329 | 0.02585 | 499.66117 | -967.98099 |
| Mixreg-Tukey | 0.48827 (0.00395) | 0.56437 (0.12794) | 0.08586 (0.05626) | 1.24795 (0.06945) | -0.08461 (0.03364) | 0.05311 | 0.02582 | 689.9801 | -1348.61901 |
| MixregGM-Mallows | 0.48932 (0.00089) | 0.56686 (0.11890) | 0.08471 (0.05394) | 1.24541 (0.10830) | -0.08274 (0.04886) | 0.04393 | 0.02451 | **772.99552** | **-1514.64969** |
| MixregGM-Schweppe | 0.48932 (0.00089) | 0.56686 (0.11890) | 0.08471 (0.05394) | 1.24541 (0.10830) | -0.08274 (0.04886) | 0.04393 | 0.02451 | 772.99553 | -1514.64939 |

Note: Value in parentheses indicates the standard errors

To see the performance of the estimators we add five outliers in the $x$ direction. The scatter plots with the fitted regression lines obtained from Mixreg-Huber, Mixreg-Tukey, MixregGM-Mallows and MixregGM-Schweppe procedures are presented in Figure 3. We also give the parameter estimates, standard errors of estimates and the ICL criterion values in Table 6. We can observe from this figure that Mixreg-Huber and Mixreg-Tukey are affected by the outliers. However, since MixregGM-Mallows and MixregGM-Schweppe are robust to the outliers in the $x$ direction, they are not influenced by the outliers. According to the ICL criterion, MixregGM-Mallows has the best fit.

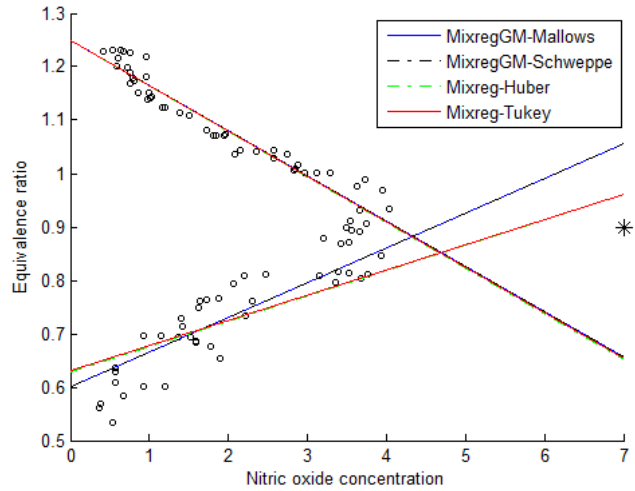

**Figure 3.** Fitted mixture regression lines for the ethanol data set with five outliers

**Table 6.** Parameter estimates, standard errors of estimates and ICL information criterion for fitting different mixture regression models to the ethanol data set with five outliers

| Model | Parameter Estimates | | | | | | | Information Criterion | |
|---|---|---|---|---|---|---|---|---|---|
| | $\hat{\pi}_1$ | $\hat{\beta}_{10}$ | $\hat{\beta}_{20}$ | $\hat{\beta}_{11}$ | $\hat{\beta}_{21}$ | $\hat{\sigma}_1$ | $\hat{\sigma}_2$ | $\ell_c(\hat{\Theta})$ | ICL |
| Mixreg-Huber | 0.48118 (0.00225) | 0.62943 (0.13828) | 0.04735 (0.04763) | 1.24886 (0.10017) | -0.08520 (0.05000) | 0.06563 | 0.02762 | 435.74227 | -839.75635 |
| Mixreg-Tukey | 0.48167 (0.00025) | 0.63033 (0.15152) | 0.04713 (0.05061) | 1.24825 (0.07439) | -0.08490 (0.03928) | 0.06551 | 0.02750 | 612.95842 | -11941.88649 |
| MixregGM-Mallows | 0.50623 (0.00372) | 0.60086 (0.13824) | 0.06515 (0.05477) | 1.24765 (0.11844) | -0.08452 (0.05746) | 0.07056 | 0.02507 | **642.35526** | **-1252.98232** |
| MixregGM-Schweppe | 0.50623 (0.00372) | 0.60086 (0.11499) | 0.06515 (0.03452) | 1.24765 (0.11844) | -0.08452 (0.05746) | 0.07056 | 0.02507 | 642.31718 | -1252.906171 |

Note: Value in parentheses indicates the standard errors.

## 5. Conclusions

In this paper, we have proposed a robust mixture regression model based on the GM regression estimation method. We have given an EM type algorithm to obtain the estimators for the proposed mixture regression model. We have provided a simulation study and a real data example to compare the performance of the proposed estimators over the estimators proposed by Bai (2010) and Bai et al. (2012). The simulation results have shown that all of the estimators behave similar in the cases of normality, heavy-tailedness without outlier in the $x$ direction. However, the mixture regression model based on the GM-estimation method outperforms the mixture regression model proposed by Bai (2010) and Bai et al. (2012) when the data includes outliers in the $x$ and $y$ directions. From the real data example, we observe the similar results.

## Appendix

Using the asymptotic covariance matrix for the mixture regression model based on the GM-estimation method, we obtain the standard errors given in Section 4. For the simplicity, we assume that scale parameters $(\sigma_1, \dots, \sigma_g)$ are fixed. Let $\widehat{\boldsymbol{\Theta}}$ be the estimate of the unknown parameter vector $\boldsymbol{\Theta} = (w_1, \cdots, w_g, \boldsymbol{\beta}_1, \cdots, \boldsymbol{\beta}_g)'$ in the mixture regression model given in (2). $\widehat{\boldsymbol{\Theta}}$ can be obtained by solving the following equations

$$\sum_{j=1}^{n} \hat{z}_{ij} x_j \eta\left(x_j, \frac{y_j - x_j' \boldsymbol{\beta}_i}{\sigma_i}\right) = 0,$$

$$\pi_i = \sum_{j=1}^{n} \frac{\hat{z}_{ij}}{n}, \quad i = 1, \dots, g,$$

where

$$\hat{z}_{ij} = \frac{\pi_i \phi(y_j; x_j' \boldsymbol{\beta}_i, \sigma_i^2)}{\sum_{i=1}^{g} \pi_i \phi(y_j; x_j' \boldsymbol{\beta}_i, \sigma_i^2)}.$$

Let $\boldsymbol{t}_j = (x_j', y_j)'$ and

$$\mathrm{H}(\boldsymbol{t}_j, \boldsymbol{\Theta}) = \left\{\hat{z}_{1j} x_j \eta\left(x_j, \frac{y_j - x_j' \boldsymbol{\beta}_1}{\sigma_1}\right), \dots, \hat{z}_{gj} x_j \eta\left(x_j, \frac{y_j - x_j' \boldsymbol{\beta}_g}{\sigma_g}\right), \hat{z}_{1j} - \pi_1, \dots, \hat{z}_{g-1,j} - \pi_{g-1}\right\}.$$

Then, the solution will be obtained using the following equation

$$S(\boldsymbol{\Theta}) = \frac{1}{n} \sum_{j=1}^{n} \mathrm{H}(\boldsymbol{t}_j, \boldsymbol{\Theta}) = 0.$$

Let

$$M = E_{\boldsymbol{\Theta}_0}\left\{\frac{\partial \mathrm{H}(\boldsymbol{t}_j, \boldsymbol{\Theta})}{\partial \boldsymbol{\Theta}'}\right\}, \tag{17}$$

$$Q = E_{\boldsymbol{\Theta}_0}\left\{\mathrm{H}(\boldsymbol{t}_j, \boldsymbol{\Theta}) \mathrm{H}(\boldsymbol{t}_j, \boldsymbol{\Theta})'\right\}. \tag{18}$$

Under the certain conditions, Maronna and Yohai (1981) show that GM-estimators are consistent and asymptotically normal. Thus, the $\widehat{\boldsymbol{\Theta}}$ has the following distribution with asymptotic covariance matrix $V$

$$\sqrt{n}(\widehat{\boldsymbol{\Theta}} - \boldsymbol{\Theta}_0) \xrightarrow{d} N(0, V),$$

where $V = M^{-1}QM^{-1}$. The standard errors of $\widehat{\boldsymbol{\Theta}}$ are obtained from the square root of the diagonal elements of the covariance matrix $V$. Note that for the finite sample case the expectation in (17) and (18) will be replaced by the average to compute the standard errors.

**References**

Azzalini, A. and Capitanio, A. 2003. Distributions generated by perturbation of symmetry with emphasis on a multivariate skew t distribution. Journal of the Royal Statistical Society, Series B (Statistical Methodology), 65(2), 367-389.

Bai, X. 2010. Robust mixture of regression models. Master Report, Kansas State University.

Bai, X., Yao, W. and Boyer, J.E. 2012. Robust fitting of mixture regression models. Computational Statistics and Data Analysis, 56(7), 2347-2359.

Bashir, S. and Carter, E.M. 2012. Robust mixture of linear regression models. Communications in Statistics-Theory and Methods, 41(18), 3371-3388.

Biernacki, C., Celeux, G. and Govaert, G. 2000. Assessing a mixture model for clustering with the integrated completed likelihood. IEEE Transactions on Pattern Analysis and Machine Intelligence, 22(3), 719-725.

Dempster, A.P., Laird, N.M. and Rubin, D.B. 1977. Maximum likelihood from incomplete data via the EM algorithm. Journal of the Royal Statistical Society, Series B (Methodological), 39, 1-38.

Doğru, F.Z. 2015. Robust parameter estimation in mixture regression models. PhD thesis, Ankara University.

Doğru, F.Z. and Arslan, O. 2015. Robust mixture regression based on the skew *t* distribution (submitted).

Doğru, F.Z. and Arslan, O. 2015. Robust mixture regression using the mixture of different distributions (revised).

Hampel, F.R. 1978. Optimally bounding the gross-error-sensitivity and the influence of position in factor space. Proceedings of the ASA Statistical Computing Section, ASA, Washington, D.C., 59-64.

Hampel, F.R., Ronchetti, E. M., Rousseeuw, P.J. and Stahel, W.A. 1986. Robust Statistics: The Approach Based on Influence Functions, Wiley, New York.

Handschin, E., Schweppe, F., Kohlas, J. and Fiechter, A. 1975. Bad data analysis for power system state estimation. IEEE Transactions on Power Apparatus and Systems, PAS-94(2), 329-337.

Huber, P.J. and Ronchetti, E.M. 2009. Robust Statistics, Wiley, New York.

Hurn, M., Justel, A. and Robert, C.P. 2003. Estimating mixtures of regressions. Journal of Computational and Graphical Statistics, 12(1), 55-79.

Hurvich, C.M., Simonoff, J.S. and Tsai, C.L. 1998. Smoothing parameter selection in nonparametric regression using an improved akaike information criterion. Journal of the Royal Statistical Society, Series B (Statistical Methodology), 60(2), 271-293.

Jurečková, J. and Picek, J. 2006. Robust statistical methods with R. Chapman & Hall/CRC, Boca Raton, FL.

Krasker, W.S. 1980. Estimation in linear regression models with disparate data points. Econometrica, 48(6), 1333-1346.

Krasker, W.S. and Welsch, R.E. 1982. Efficient bounded-influence regression estimation. Journal of the American Statistical Association, 77(379), 595-604.

Lin, T.I., Lee, J.C. and Hsieh, W.J. 2007a. Robust mixture modeling using the skew t distribution. Statistics and Computing, 17, 81-92.

Loader, C.R. 1999. Local Regression and Likelihood. Statistics and Computing Series. Springer Verlag, New York.

Mallows, C.L. 1975. On some topics in robustness. Technical memorandum, Bell Telephone Laboratories, Murray Hill, N.J.

Markatou, M. 2000. Mixture models, robustness, and the weighted likelihood methodology. Biometrics, 56(2), 483-486.

Maronna, R.A., Martin, R.D. and Yohai, V.J. 2006. Robust Statistics: Theory and Methods. Wiley, New York.

Maronna, R.A. and Yohai, V.J. 1981. Asymptotic behavior of general M-estimates for regression and scale with random carriers. Zeitschrift für Wahrscheinlichkeitstheorie und Verwandte Gebiete, 58(1), 7-20.

Neykov, N., Filzmoser, P., Dimova, R. and Neytchev, P. 2007. Robust fitting of mixtures using the trimmed likelihood estimator. Computational Statistics and Data Analysis, 52(1), 299-308.

Quandt, R.E. 1972. A new approach to estimating switching regressions. Journal of the American Statistical Association, 67(338), 306-310.

Quandt, R.E. and Ramsey, J.B. 1978. Estimating mixtures of normal distributions and switching regressions. Journal of the American Statistical Association, 73(364), 730-752.

Rousseeuw, P.J. and Van Driessen, K. 1999. A fast algorithm for the minimum covariance determinant estimator. Technometrics, 41(3), 212-223.

Shen, H., Yang, J. and Wang, S. 2004. Outlier detecting in fuzzy switching regression models. Artificial Intelligence: Methodology, Systems, and Applications. In Lecture Notes in Computer Science, 3192, 208-215.

Simpson, D.G., Ruppert, D. and Carroll, R.J. 1992. On one-step GM estimates and stability of inferences in linear regression. Journal of the American Statistical Association, 87(418), 439-450.

Song, W., Yao, W. and Xing, Y. 2014. Robust mixture regression model fitting by Laplace distribution. Computational Statistics and Data analysis, 71, 128-137.

Wei, Y. 2012. Robust mixture regression models using t-distribution. Master Report, Kansas State University.

Yao, W., Wei, Y. and Yu, C. 2014. Robust mixture regression using the t-distribution. Computational Statistics and Data Analysis, 71, 116–127.

Zhang, J. 2013. Robust mixture regression modeling with Pearson Type VII distribution. Master Report, Kansas State University.